\documentclass[12pt]{article}%
\usepackage{amsfonts}
\usepackage{sw20bams}
\usepackage{amsmath}
\usepackage{amssymb}
\usepackage{graphicx}%
\providecommand{\U}[1]{\protect\rule{.1in}{.1in}}
\begin{document}

\title{Conjectures about Traffic Light Queues}
\author{Steven Finch and Guy Louchard}
\date{October 31, 2018}
\maketitle

\begin{abstract}
In discrete time, $\ell$-blocks of red lights are separated by $\ell$-blocks
of green lights. Cars arrive at random. \ The maximum line length of idle cars
is fully understood for $\ell=1$, but only partially for $2\leq\ell\leq3$.

\end{abstract}

\footnotetext{Copyright \copyright \ 2018 by Steven R. Finch. All rights
reserved.}Let $\ell\geq1$ be an integer. \ Let $X_{0}=0$ and $X_{1}$, $X_{2}$,
\ldots, $X_{n}$ be a sequence of independent random variables satisfying%
\[%
\begin{array}
[c]{ccccc}%
\mathbb{P}\left\{  X_{i}=1\right\}  =p, &  & \mathbb{P}\left\{  X_{i}%
=0\right\}  =q &  & \text{if }i\equiv1,2,\ldots,\ell\,\operatorname{mod}%
\,2\ell;
\end{array}
\]%
\[%
\begin{array}
[c]{ccccc}%
\mathbb{P}\left\{  X_{i}=0\right\}  =p, &  & \mathbb{P}\left\{  X_{i}%
=-1\right\}  =q &  & \text{if }i\equiv\ell+1,\ell+2,\ldots,2\ell
\,\operatorname{mod}\,2\ell
\end{array}
\]
for each $1\leq i\leq n$. \ Define $S_{0}=X_{0}$ and $S_{j}=\max\left\{
S_{j-1}+X_{j},0\right\}  $ for all $1\leq j\leq n$. \ Thus cars arrive at a
one-way intersection according to a Bernoulli($p$) distribution; \ when the
signal is red ($1\leq i\leq\ell$), no cars may leave; \ when the signal is
green ($\ell+1\leq i\leq2\ell$), a car must leave (if there is one). The
quantity $M_{n}=\max\nolimits_{0\leq j\leq n}S_{j}$ is the worst-case traffic
congestion (as opposed to the average-case often cited). \ Only the
circumstance when $\ell=1$ is amenable to rigorous treatment, as far as is
known. Let%
\[
\chi_{1}(p)=\frac{p(q-p)^{2}}{q^{3}}%
\]
\ where $q=1-p$. \ It is proved in \cite{Fi1-tlqc} via a theorem in
\cite{HL-tlqc, LP-tlqc} that%
\[
\mathbb{P}\left\{  M_{n}\leq\log_{q^{2}/p^{2}}(n)+h\right\}  \sim\exp\left[
-\frac{\chi_{1}(p)}{2}\left(  \frac{q^{2}}{p^{2}}\right)  ^{-h}\right]  ,
\]%
\[%
\begin{array}
[c]{ccc}%
\mathbb{E}\left(  M_{n}\right)  \sim\dfrac{\ln(n)}{\ln\left(  \frac{q^{2}%
}{p^{2}}\right)  }+\dfrac{\gamma+\ln\left(  \frac{\chi_{1}(p)}{2}\right)
}{\ln\left(  \frac{q^{2}}{p^{2}}\right)  }+\dfrac{1}{2}+\varphi(n), &  &
\mathbb{V}\left(  M_{n}\right)  \sim\dfrac{\pi^{2}}{6}\dfrac{1}{\ln\left(
\frac{q^{2}}{p^{2}}\right)  ^{2}}+\dfrac{1}{12}+\psi(n)
\end{array}
\]
as $n\rightarrow\infty$, assuming $p<q$. \ The symbol $\gamma$ denotes Euler's
constant \cite{Fi2-tlqc} ; $\varphi$ and $\psi$ are periodic functions of
$\log_{q^{2}/p^{2}}(n)$ with period $1$ and of small amplitude.

For $\ell=2$, we have the following conjecture:%
\[
\mathbb{P}\left\{  M_{n}\leq\log_{q^{2}/p^{2}}(n)+h\right\}  \sim\exp\left[
-\frac{\chi_{2}(p)}{4}\left(  \frac{q^{2}}{p^{2}}\right)  ^{-h}\right]
\]
where%
\[
\chi_{2}(p)=\frac{(q-p)^{2}}{4q^{6}}\left[  \left(  1-8p^{2}+16p^{3}%
-8p^{4}\right)  +(q-p)\sqrt{1+4pq}\right]  .
\]
For $\ell=3$, we likewise have:%
\[
\mathbb{P}\left\{  M_{n}\leq\log_{q^{2}/p^{2}}(n)+h\right\}  \sim\exp\left[
-\frac{\chi_{3}(p)}{6}\left(  \frac{q^{2}}{p^{2}}\right)  ^{-h}\right]
\]
where%
\[
\chi_{3}(p)=\frac{(q-p)^{2}}{12pq^{9}}\left[  a+(q-p)^{2}b\theta+(q-p)\sqrt
{2}\sqrt{c+ab\theta}\right]  ,
\]%
\[
a=1-4p+10p^{2}-52p^{3}+226p^{4}-520p^{5}+640p^{6}-400p^{7}+100p^{8},
\]%
\[
b=1-2p+6p^{2}-8p^{3}+4p^{4},
\]%
\begin{align*}
c  &  =1-4p+16p^{2}-104p^{3}+506p^{4}-1808p^{5}+5604p^{6}-15576p^{7}%
+35574p^{8}\\
&  -61160p^{9}+75152p^{10}-63440p^{11}+34840p^{12}-11200p^{13}+1600p^{14},
\end{align*}%
\[
\theta=\sqrt{1+4pq+16p^{2}q^{2}}.
\]
Our ad hoc technique for deriving such formulas is based on the numerical
solution of a large determinantal equation, followed by the recognition of
algebraic quantities given high-precision decimals. \ A better-justified
method is currently being developed \cite{LF-tlqc}.

\section{Computational Technique}

Let $k\geq1$ be an integer. \ Define $(k+1)\times(k+1)$ matrices%
\[
U_{k}=\left(
\begin{array}
[c]{cccccccc}%
q & p & 0 & 0 & \cdots & 0 & 0 & 0\\
0 & q & p & 0 & \cdots & 0 & 0 & 0\\
0 & 0 & q & p & \cdots & 0 & 0 & 0\\
0 & 0 & 0 & q & \cdots & 0 & 0 & 0\\
\vdots & \vdots & \vdots & \vdots & \ddots & \vdots & \vdots & \vdots\\
0 & 0 & 0 & 0 & \cdots & q & p & 0\\
0 & 0 & 0 & 0 & \cdots & 0 & q & p\\
0 & 0 & 0 & 0 & \cdots & 0 & 0 & q
\end{array}
\right)  ,
\]%
\[
V_{k}=\left(
\begin{array}
[c]{cccccccc}%
1 & 0 & 0 & 0 & \cdots & 0 & 0 & 0\\
q & p & 0 & 0 & \cdots & 0 & 0 & 0\\
0 & q & p & 0 & \cdots & 0 & 0 & 0\\
0 & 0 & q & p & \cdots & 0 & 0 & 0\\
\vdots & \vdots & \vdots & \vdots & \ddots & \vdots & \vdots & \vdots\\
0 & 0 & 0 & 0 & \cdots & p & 0 & 0\\
0 & 0 & 0 & 0 & \cdots & q & p & 0\\
0 & 0 & 0 & 0 & \cdots & 0 & q & p
\end{array}
\right)
\]
and let $z_{k}$ denote the (positive real) solution of the equation%
\[
\det\left[  I-U_{k}^{\ell}V_{k}^{\ell}z\right]  =0
\]
that is closest to unity. \ The quantity $z_{k}$ can be numerically estimated.
\ Define%
\[
\chi_{\ell}(p)=\lim_{k\rightarrow\infty}\frac{z_{k}-1}{(p/q)^{2k}}.
\]
For suitably large $k$ and sufficiently accurate $z_{k}$, the preceding ratio
can be recognized as a specific algebraic number when $p$ is rational. \ From
a set of such algebraic numbers and corresponding $p$ values, a $\chi_{\ell
}(p)$ formula can be inferred.

For example, when $\ell=3$ and $p=1/3$, taking $k\approx400$ and employing
$100$ precise digits gives
\[
\chi_{\ell}(1/3)=\frac{1393+61\sqrt{217}+\sqrt{2416130+169946\sqrt{217}}%
}{6144};
\]
when instead $p=1/5$, we obtain%
\begin{align*}
\chi_{\ell}(1/5)  &  =\frac{27\left[  18025+489\sqrt{1281}+5\sqrt
{25206642+705138\sqrt{1281}}\right]  }{1048576}\\
&  =\frac{9\left[  162225+4401\sqrt{1281}+3\sqrt{5671494450+158656050\sqrt
{1281}}\right]  }{3145728}.
\end{align*}
In contrast, when $p=1/17$, taking $k\approx2500$ and employing $300$ precise
digits gives%

\begin{align*}
\chi_{\ell}(1/17)  &  =\tfrac{675\left[  613160569+1882425\sqrt{106113}%
+5\sqrt{30079190568067506+92338302727986\sqrt{106113}}\right]  }%
{274877906944}\\
&  =\tfrac{225\left[  5518445121+16941825\sqrt{106113}+15\sqrt
{270712715112607554+831044724551874\sqrt{106113}}\right]  }{824633720832};
\end{align*}
when instead $p=1/19$, we obtain%
\[
\chi_{\ell}(1/19)=\tfrac{289\left[  13775887153+34281469\sqrt{161497}%
+17\sqrt{1313388976733016770+3268219019952026\sqrt{161497}}\right]
}{2380311484416}.
\]

It is not hard to surmise from $xy$-data (with $x$ equal to $1/p$ odd)
\[
\left(
\begin{array}
[c]{c}%
3\\
6144
\end{array}
\right)  ,\left(
\begin{array}
[c]{c}%
5\\
3145728
\end{array}
\right)  ,\ldots,\left(
\begin{array}
[c]{c}%
17\\
824633720832
\end{array}
\right)  ,\left(
\begin{array}
[c]{c}%
19\\
2380311484416
\end{array}
\right)  ,\ldots
\]
that the denominator is prescribed by%
\[
12(-1+x)^{9}.
\]
The coefficient $a$ is found via polynomial regression on $xy$-data%
\[
\left(
\begin{array}
[c]{c}%
3\\
1393
\end{array}
\right)  ,\left(
\begin{array}
[c]{c}%
5\\
162225
\end{array}
\right)  ,\ldots,\left(
\begin{array}
[c]{c}%
17\\
5518445121
\end{array}
\right)  ,\left(
\begin{array}
[c]{c}%
19\\
13775887153
\end{array}
\right)  ,\ldots
\]
yielding%
\[
y=100-400x+640x^{2}-520x^{3}+226x^{4}-52x^{5}+10x^{6}-4x^{7}+x^{8}.
\]
The coefficient $b$ is found via data%
\[
\left(
\begin{array}
[c]{c}%
3\\
61
\end{array}
\right)  ,\left(
\begin{array}
[c]{c}%
5\\
4401
\end{array}
\right)  ,\ldots,\left(
\begin{array}
[c]{c}%
17\\
16941825
\end{array}
\right)  ,\left(
\begin{array}
[c]{c}%
19\\
34281469
\end{array}
\right)  ,\ldots
\]
yielding%
\begin{align*}
y  &  =16-48x+60x^{2}-40x^{3}+18x^{4}-6x^{5}+x^{6}\\
&  =(-2+x)^{2}(4-8x+6x^{2}-2x^{3}+x^{4}).
\end{align*}
The coefficient $c$ is found via data%
\[
\left(  {\scriptsize
\begin{array}
[c]{c}%
3\\
2416130
\end{array}
}\right)  ,\left(  {\scriptsize
\begin{array}
[c]{c}%
5\\
5671494450
\end{array}
}\right)  ,\ldots,\left(  {\scriptsize
\begin{array}
[c]{c}%
17\\
270712715112607554
\end{array}
}\right)  ,\left(  {\scriptsize
\begin{array}
[c]{c}%
19\\
1313388976733016770
\end{array}
}\right)  ,\ldots
\]
yielding%
\begin{align*}
c  &  =3200-22400x+69680x^{2}-126880x^{3}+150304x^{4}-122320x^{5}+71148x^{6}\\
&  -31152x^{7}+11208x^{8}-3616x^{9}+1012x^{10}-208x^{11}+32x^{12}%
-8x^{13}+2x^{14}.
\end{align*}
Little insight is provided by a brute-force technique, but effectiveness is
key. \ The greatest difficulty is knowing when cancellation of common factors
has occurred and hence needs remedy (as when $p=1/5$ and $p=1/17$). \ Semi-log
plots for the coefficients were helpful in detecting such. \ This issue is
more challenging still when $1/p$ is even, e.g., when $p=1/10$:%

\begin{align*}
\chi_{\ell}(1/10)  &  =\tfrac{64\left[  16650025+3818752\sqrt{19}%
+80\sqrt{86608486817+19869473834\sqrt{19}}\right]  }{1162261467}\\
&  =\tfrac{64\left[  66600100+545536\sqrt{14896}+8\sqrt
{138573578907200+1135398504800\sqrt{14896}}\right]  }{4649045868}.
\end{align*}
Both a factor of $28$ is transferred under the radical (of $19$) and a factor
of $4$ is introduced throughout. This similarly happens when $p=1/12$:%

\begin{align*}
\chi_{\ell}(1/12) &  =\tfrac{100\left[  76862569+12636400\sqrt{37}%
+20\sqrt{29539863834326+4856330834558\sqrt{37}}\right]  }{7073843073}\\
&  =\tfrac{100\left[  307450276+1805200\sqrt{29008}+10\sqrt
{1890551285396864+11100184764704\sqrt{29008}}\right]  }{28295372292}.
\end{align*}
Of course, the formulas require checking for non-integer $1/p$. \ The case
$\ell=2$ is less tediously studied and will additionally be the focus of
\cite{LF-tlqc}.

\section{Queue Data}

Let $n=10^{10}$. \ For each $p\in\{1/5,1/3\}$, we\ generated $40000$ traffic
light queues (for both $\ell=2$ and $\ell=3$) and produced an empirical
histogram for the maximum $M_{n}$. \ Figures 1--4 contain these histograms (in
blue) along with our theoretical predictions (in red). \ The fit is excellent.
\ Identical barcharts appeared in \cite{Fi1-tlqc} but with ill-informed
predictions emerging from $\ell=1$. \
\begin{figure}[ptb]%
\centering
\includegraphics[
height=3.0251in,
width=3.691in
]%
{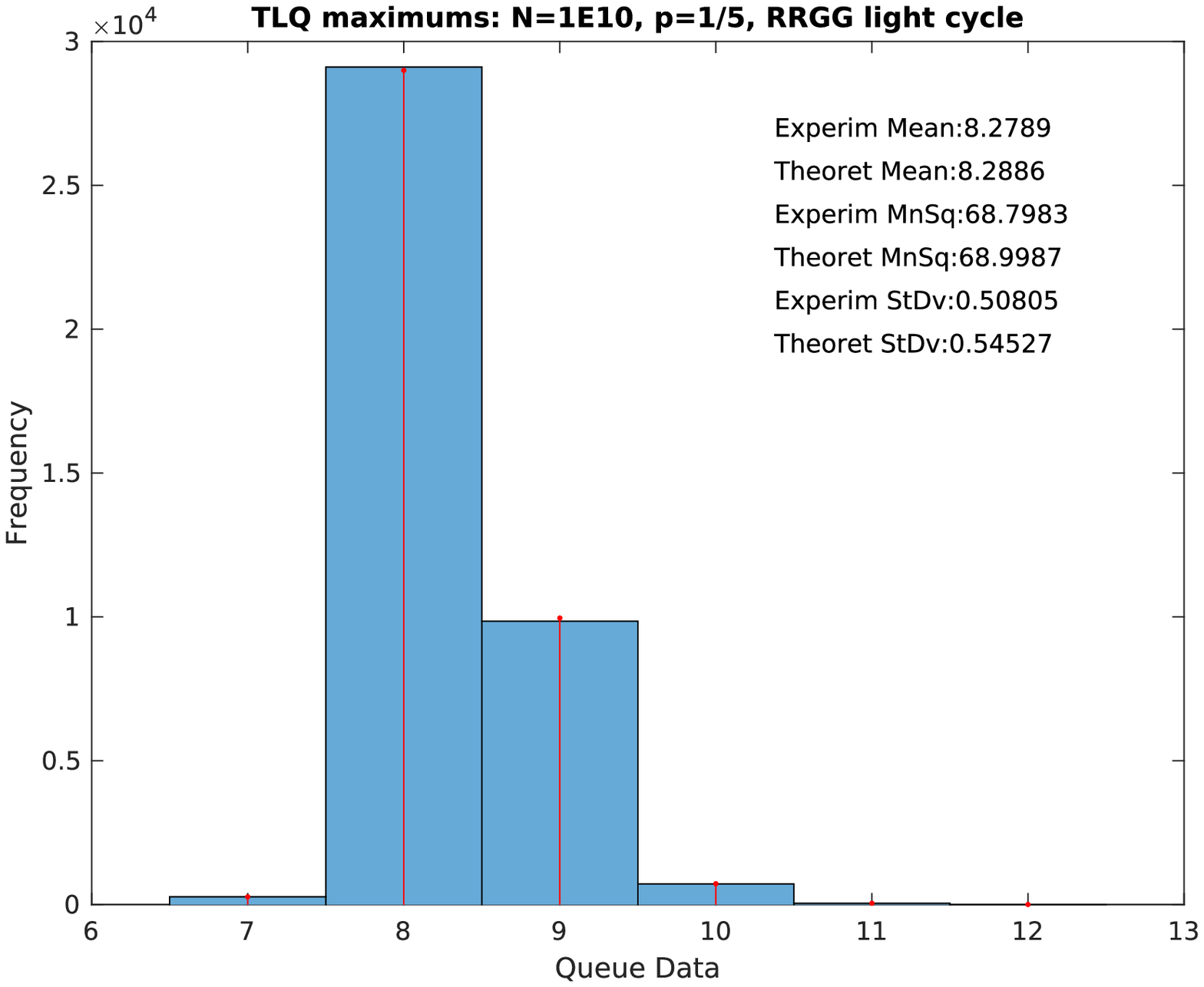}%
\caption{ }%
\end{figure}
\begin{figure}[ptb]%
\centering
\includegraphics[
height=3.0251in,
width=3.8156in
]%
{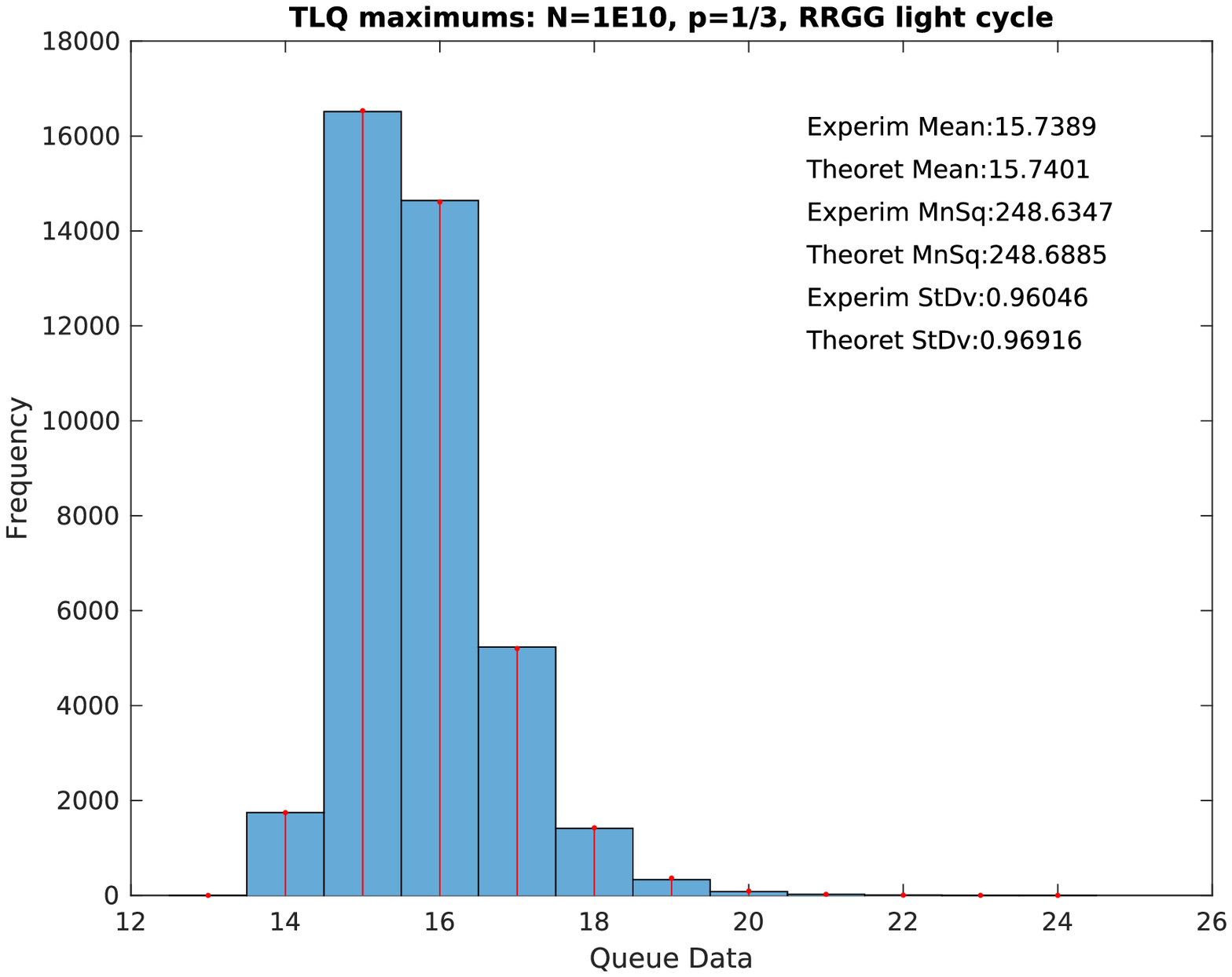}%
\caption{ \ }%
\end{figure}
\begin{figure}[ptb]%
\centering
\includegraphics[
height=3.0251in,
width=3.691in
]%
{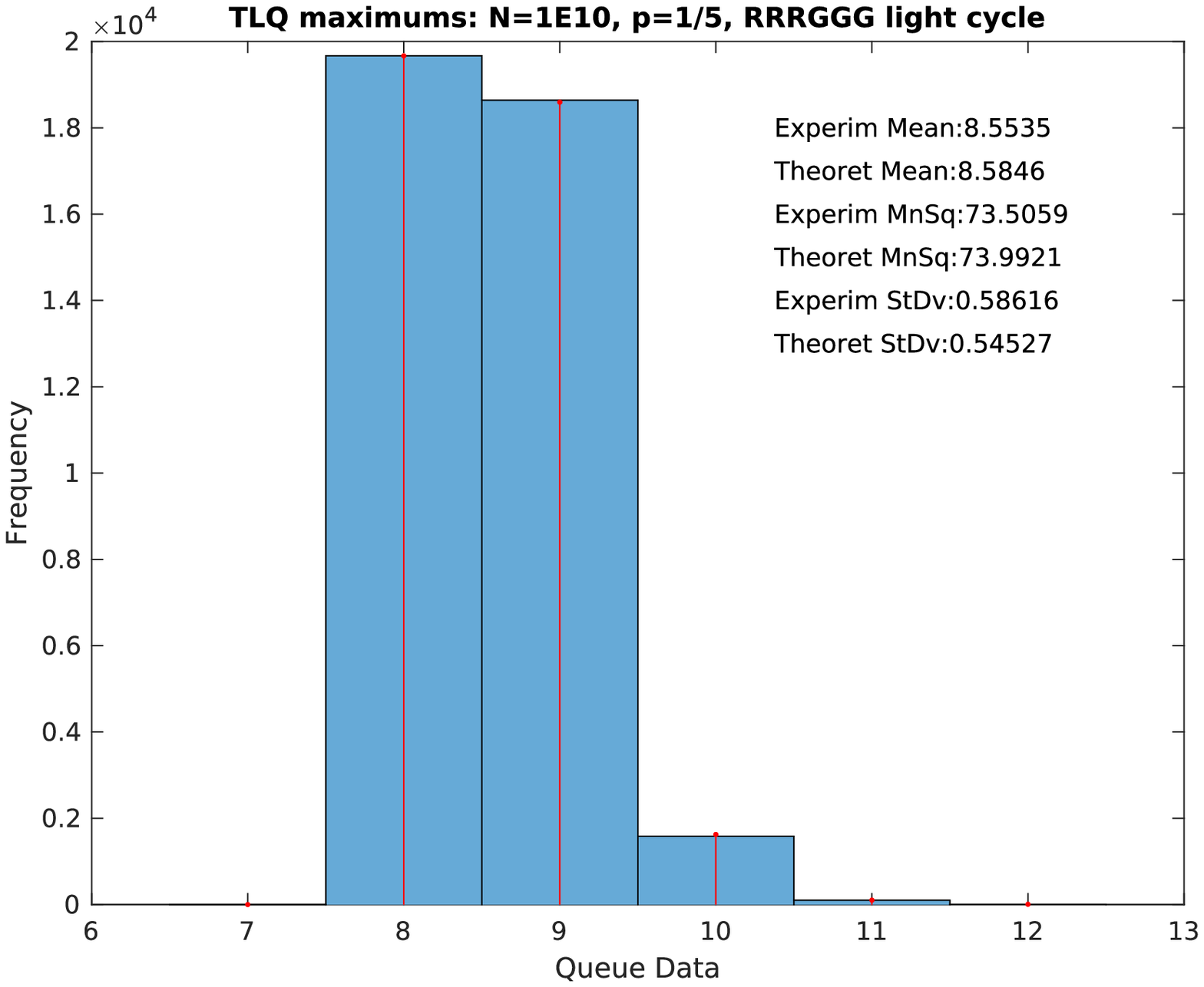}%
\caption{ }%
\end{figure}
\begin{figure}[ptb]%
\centering
\includegraphics[
height=3.0251in,
width=3.8156in
]%
{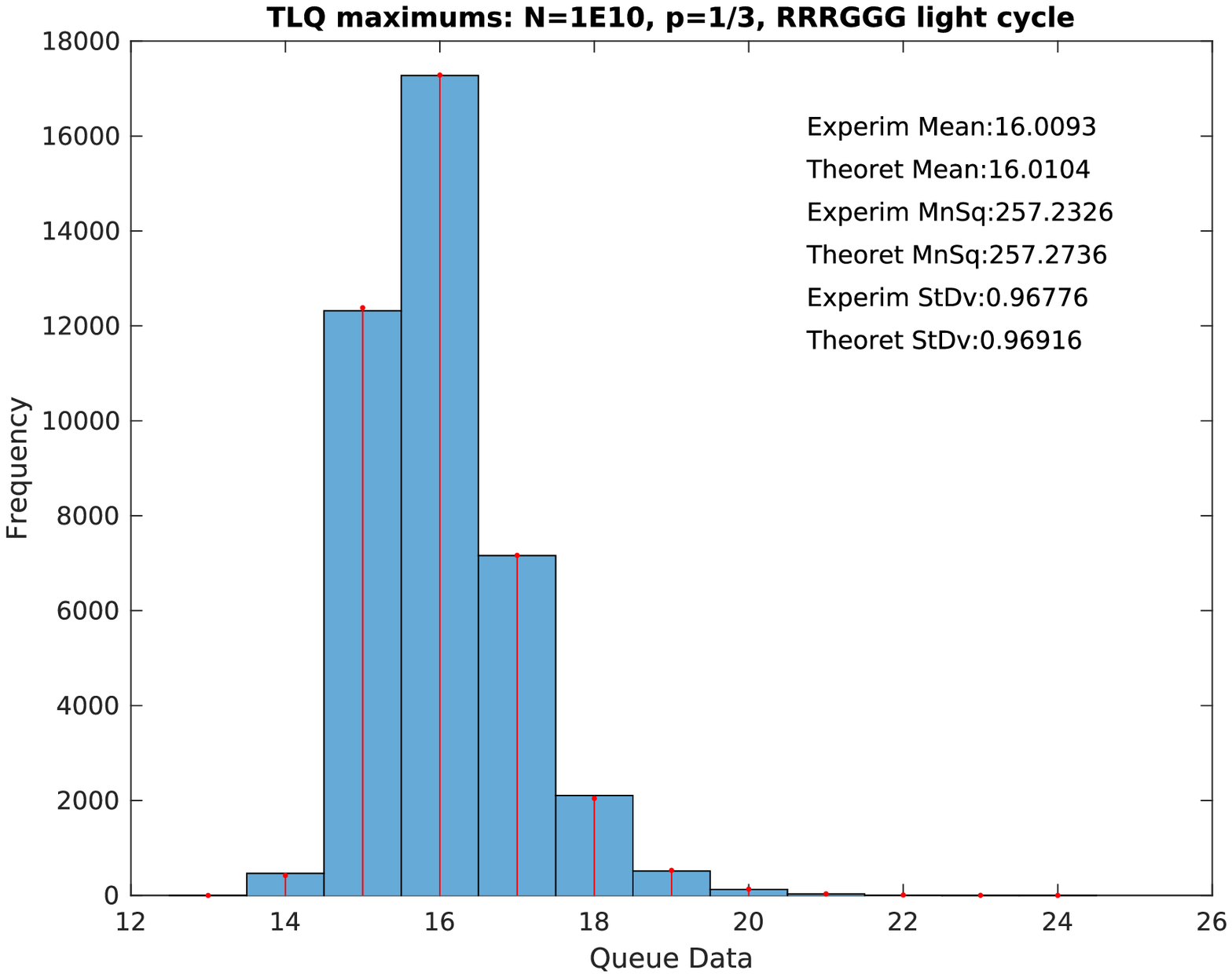}%
\caption{ }%
\end{figure}

\section{Random Lights}

Let $Y_{1}$, $Y_{2}$, \ldots, $Y_{n}$ be a sequence of independent
Bernoulli$(1/2)$ variables. \ Rather than defining increments $X_{i}$
deterministically based on $i\operatorname{mod}\,2\ell$, let us define $X_{i}$
randomly based on $Y_{i}$ as follows:%
\[%
\begin{array}
[c]{ccccc}%
\mathbb{P}\left\{  X_{i}=1\right\}  =p, &  & \mathbb{P}\left\{  X_{i}%
=0\right\}  =q &  & \text{if }Y_{i}=1;
\end{array}
\]%
\[%
\begin{array}
[c]{ccccc}%
\mathbb{P}\left\{  X_{i}=0\right\}  =p, &  & \mathbb{P}\left\{  X_{i}%
=-1\right\}  =q &  & \text{if }Y_{i}=0\text{.}%
\end{array}
\]
The corresponding sequence $S_{1}$, $S_{2}$, \ldots, $S_{n}$, reflected at the
origin as before, is a lazy random walk with expected maximum approximately
equal to \cite{Fi3-tlqc}%
\[
E_{0}\left(  n,p\right)  =\dfrac{\ln\left(  \frac{n}{2}\right)  }{\ln\left(
\frac{q}{p}\right)  }+\dfrac{\gamma+\ln\left(  \frac{\chi_{0}(p)}{1}\right)
}{\ln\left(  \frac{q}{p}\right)  }+\dfrac{1}{2}%
\]
where%
\[
\chi_{0}(p)=\frac{p(q-p)^{2}}{q^{2}}.
\]
Define also expected maximums associated with $\ell=1$, $2$, $3$:
\[
E_{1}\left(  n,p\right)  =\dfrac{\ln(n)}{\ln\left(  \frac{q^{2}}{p^{2}%
}\right)  }+\dfrac{\gamma+\ln\left(  \frac{\chi_{1}(p)}{2}\right)  }%
{\ln\left(  \frac{q^{2}}{p^{2}}\right)  }+\dfrac{1}{2},
\]%
\[
E_{2}\left(  n,p\right)  =\dfrac{\ln(n)}{\ln\left(  \frac{q^{2}}{p^{2}%
}\right)  }+\dfrac{\gamma+\ln\left(  \frac{\chi_{2}(p)}{4}\right)  }%
{\ln\left(  \frac{q^{2}}{p^{2}}\right)  }+\dfrac{1}{2},
\]%
\[
E_{3}\left(  n,p\right)  =\dfrac{\ln(n)}{\ln\left(  \frac{q^{2}}{p^{2}%
}\right)  }+\dfrac{\gamma+\ln\left(  \frac{\chi_{3}(p)}{6}\right)  }%
{\ln\left(  \frac{q^{2}}{p^{2}}\right)  }+\dfrac{1}{2}.
\]
Fix $n=10^{10}$ for sake of definiteness. \ Figure 5 shows that $\ell=1$
possesses the best $E_{\ell}\left(  n,p\right)  $, in the sense of minimizing
traffic backup, however only slightly compared against $\ell=2$ and $\ell=3$.
\ It may be surprising that $\ell=0$ possesses the worst $E_{\ell}\left(
n,p\right)  $ by far. \ Is there a red/green light strategy (either
deterministic or random) that both maintains equiprobability and yet improves
upon $\ell=1$? \ An answer to this question would be good to see someday.%
\begin{figure}[ptb]%
\centering
\includegraphics[
height=2.8582in,
width=5.3039in
]%
{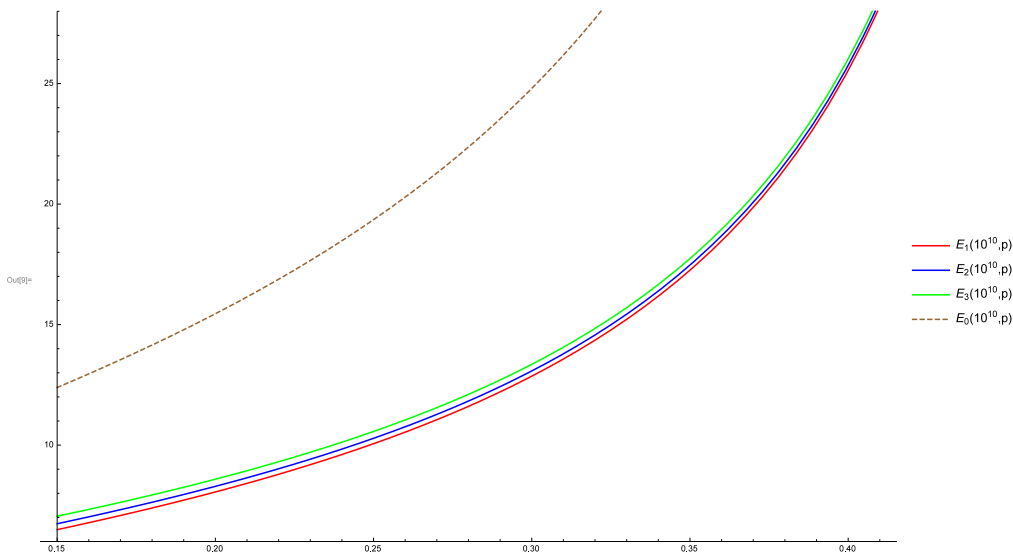}%
\caption{ Expected maximums (approximate) as functions of $0.15<p<0.41$.}%
\end{figure}

\section{Acknowledgements}

The creators of Julia, Mathematica and Matlab, as well as administrators of
the MIT\ Engaging Cluster, earn our gratitude every day.


\begin{thebibliography}{9}                                                                                                %


\bibitem {Fi1-tlqc}S. Finch, The maximum of an asymmetric simple random walk
with reflection, arXiv:1808.01830.

\bibitem {HL-tlqc}P. Hitczenko and G. Louchard, Distinctness of compositions
of an integer: a probabilistic analysis, \textit{Random Structures Algorithms}
19 (2001) 407--437; MR1871561.

\bibitem {LP-tlqc}G. Louchard and H. Prodinger, Asymptotics of the moments of
extreme-value related distribution functions, \textit{Algorithmica} 46 (2006)
431--467 (long version available at http://www.ulb.ac.be/di/mcs/louchard/); MR2291964.

\bibitem {Fi2-tlqc}S. Finch, Euler-Mascheroni constant, \textit{Mathematical
Constants}, Cambridge Univ. Press, 2003, pp. 28--40; MR2003519.

\bibitem {LF-tlqc}S. Finch and G. Louchard, Traffic light queues and the
Poisson clumping heuristic, arXiv:1810.12058.

\bibitem {Fi3-tlqc}S. Finch, How far might we walk at random?, arXiv:1802.04615.%

\begin{tabular}
[c]{llllll}
& Steven Finch &  &  & Guy Louchard & \\
& MIT Sloan School of Management &  &  & Universit\'{e} Libre de Bruxelles &
\\
& Cambridge, MA, USA &  &  & Bruxelles, Belgium & \\
& \textit{steven\_finch@harvard.edu} &  &  & \textit{louchard@ulb.ac.be} &
\end{tabular}

\end{thebibliography}
\end{document}